\documentclass[12pt]{article}
\usepackage{amssymb}
\usepackage{amsmath,amsthm}
\usepackage[latin1]{inputenc}
\usepackage{color}
\usepackage{graphicx}

\usepackage{hyperref}

\hypersetup{colorlinks=true, urlcolor=blue}

\usepackage{epsfig}

 \newtheorem{remark}{Remark}

 \newtheorem{theorem}[remark]{Theorem}
 
 \newtheorem{corollary}[remark]{Corollary}

\title{Defensive $k$-alliances in graphs}

\author{Juan A. Rodr\'{\i}guez-Vel\'{a}zquez\footnote{juanalberto.rodriguez\@@urv.cat} and Ismael
G. Yero\footnote{ismael.gonzalez\@@urv.cat} \\
{\em Department of Computer Engineering and Mathematics}\\
Rovira i Virgili University of Tarragona\\ Av. Pa\"{\i}sos Catalans
26, 43007 Tarragona, Spain\\
Jos\'{e} M. Sigarreta\footnote{josemaria.sigarreta\@@uc3m.es}\\
{\em Departamento de Matem\'{a}ticas}\\
Universidad Carlos III de Madrid\\ Avda. de la Universidad 30, 28911
Leganés (Madrid),  Spain }

\date{}

\begin{document}

\maketitle

\begin{abstract}
Let   $G=(V,E)$ be a simple graph of order $n$ and  degree sequence
$\delta_{1}\geq \delta_{2}\geq \cdots \geq\delta_{n}$. For a
nonempty set $X\subseteq V$, and a vertex $v\in V$, $\delta_{X}(v)$
denotes the number of
  neighbors  $v$ has in $X$.
A nonempty set $S\subseteq V$ is a \emph{defensive  $k$-alliance} in
$G=(V,E)$ if $\delta _S(v)\ge \delta_{\bar{S}}(v)+k,$ $\forall v\in
S.$
 The defensive $k$-alliance number
of $G$, denoted by $a_k(G)$, is defined as the minimum cardinality
of a defensive $k$-alliance in $G$. We study the mathematical
properties of $a_k(G)$. We show that
$\left\lceil\frac{\delta_n+k+2}{2}\right\rceil \le a_k(G)\leq
n-\left\lfloor\frac{\delta_n-k}{2}\right\rfloor$ and $ a_k(G)\ge
\left\lceil\frac{n(\mu+k+1)}{n+\mu}\right\rceil,$ where $\mu$ is the
algebraic connectivity of $G$ and $k\in
\{-\delta_n,\dots,\delta_1\}$. Moreover, we show that for every
$k,r\in \mathbb{Z}$ such that $ -\delta_n\le k\le \delta_1$
 and $0 \le r \le \frac{k+\delta_n}{2}$, $
a_{_{k-2r}}(G)+r\leq a_{_k}(G)$ and, as a consequence, we show that
for every $k\in \{-\delta_n,\dots,0\}$,
 $a_k(G)\leq \left\lceil\displaystyle\frac{n+k+1}{2}\right\rceil
 .$  In the case of the line graph
${\cal L}(G)$ of a simple graph $G$, we obtain bounds on $a_k({\cal
L}(G))$ and, as a consequence of the study, we show that for any
$\delta$-regular graph, $\delta>0$, and for every $k\in
\{2(1-\delta),\dots, 0\}$, $a_k({\cal
L}(G))=\delta+\left\lceil\frac{k}{2}\right\rceil .$ Moreover, for
any ($\delta_1,\delta_2$)-semiregular bipartite graph $G$,
$\delta_1>\delta_2$, and for every $k\in \{2-\delta_1-\delta_2,\dots
,\delta_1-\delta_2\}$, $a_k({\cal L}(G))=\left\lceil
\frac{\delta_1+\delta_2+k}{2} \right\rceil.$
\end{abstract}



\section{Introduction}

The mathematical properties of alliances in graphs were first
studied by P. Kristiansen, S. M. Hedetniemi and S. T. Hedetniemi
\cite{alliancesOne}. They proposed diffe\-rent types of alliances:
namely, defensive alliances
\cite{note,GlobalalliancesOne,alliancesOne,albesiga8}, offensive
alliances \cite{fava,albesiga6,albesiga7} and dual alliances or
po\-werful alliances \cite{cota}. A generalization of these
alliances called $k$-alliances was presented by K. H. Shafique and
R. D. Dutton \cite{kdaf,kdaf1}.

In this paper, we  study the mathematical properties of defensive
$k$-alliances. We begin by stating the  terminology used. Throughout
this article, $G=(V,E)$ denotes a simple graph of order $|V|=n$ and
size $|E|=m$. We denote two adjacent vertices $u$ and $v$ by $u\sim
v$. For a nonempty set $X\subseteq V$, and a vertex $v\in V$,
 $N_X(v)$ denotes the set of neighbors  $v$ has in $X$:
$N_X(v):=\{u\in X: u\sim v\}$ and the degree of $v$ in $ X$ will be
denoted by $\delta_{X}(v)=|N_{X}(v)|.$ We denote the degree of a
vertex $v_i\in V$  by $\delta(v_i)$ (or by $\delta_i$ for short) and
the degree sequence of
 $G$ by $\delta_{1}\geq
\delta_{2}\geq \cdots \geq\delta_{n}$. The subgraph induced by
$S\subset V$ will be denoted by  $\langle S\rangle $ and the
complement of the set $S$ in $V$ will be denoted by $\bar{S}$.

A nonempty set $S\subseteq V$ is a \emph{defensive  $k$-alliance} in
$G=(V,E)$,  $k\in \{-\delta_1,\dots,\delta_1\}$, if for every $ v\in
S$,
\begin{equation}\label{cond-A-Defensiva} \delta _S(v)\ge \delta_{\bar{S}}(v)+k.\end{equation}

A vertex $v\in S$ is said to be $k$-\emph{satisfied} by the set $S$
if   (\ref{cond-A-Defensiva}) holds. Notice that
(\ref{cond-A-Defensiva}) is equivalent to
\begin{equation}\label{cond-A-Defensiva1} \delta (v)\ge 2\delta_{\bar{S}}(v)+k.\end{equation}


A defensive $(-1)$-alliance is a \emph{defensive alliance} and a
defensive $0$-alliance is a \emph{strong defensive alliance} as
defined in \cite{alliancesOne}. A defensive $0$-alliance is also
known as a \emph{cohesive set} \cite{porsi}.

Defensive alliances are the mathematical model of  web communities.
Adopting the definition of Web community proposed recently by Flake,
Law\-rence, and Giles \cite{flake}, ``a Web community is a set of
web pages having more hyperlinks (in either direction) to members of
the set than to non-members".

\section{Defensive $k$-alliance number}

The \emph{defensive $k$-alliance number} of $G$, denoted by
$a_k(G)$, is defined as the minimum cardinality of a defensive
$k$-alliance in $G$. Notice that
\begin{equation}a_{k+1}(G)\ge a_k(G).\end{equation}

The defensive $(-1)$-alliance number of $G$ is known as the
\emph{alliance number} of $G$
 and the defensive  $0$-alliance number is known as  the \emph{strong alliance number},
  \cite{alliancesOne, note, GlobalalliancesOne}.
For instance, in the case of the $3$-cube graph, $G=Q_{3}$, every\-
set composed by two adjacent vertices is a defensive alliance of
minimum cardinality and every\-  set composed by four vertices whose
induced subgraph is isomorphic to the cycle  $C_4$ is a strong
defensive alliance of minimum cardinality. Thus, $a_{-1}(Q_3)=2$ and
$a_{0}(Q_3)=4$.

 Notice that if every vertex of $G$ has even
 degree and $k$ is odd, $k=2l-1$ ,
 then
every defensive $(2l-1)$-alliance in $G$ is a defensive
$(2l)$-alliance. Hence, in such a case, $a_{2l-1}(G)=a_{2l}(G)$.
Analogously, if every vertex of $G$ has odd
 degree and $k$ is even, $k=2l$,
 then
every defensive $(2l)$-alliance in $G$ is a defensive
$(2l+1)$-alliance. Hence, in such a case, $a_{2l}(G)=a_{2l+1}(G)$.

For some graphs, there are some values of  $k\in
\{-\delta_1,\dots,\delta_1\}$,   such that defensive $k$-alliances
do not exist. For instance, for $k\ge 2$ in the case of the star
graph $S_n$, defensive $k$-alliances do not exist. By
(\ref{cond-A-Defensiva1}) we conclude that, in any graph, there are
defensive $k$-alliances for all $k\in \{-\delta_1,\dots,\delta_n\}$.
For instance, a defensive $(\delta_n)$-alliance in $G=(V,E)$ is $V$.
Moreover, if $v\in V$ is a vertex of minimum degree,
$\delta(v)=\delta_n$, then $S=\{v\}$ is a defensive $k$-alliance for
every $k\le -\delta_n$. As $a_k(G)=1$ for $k\le -\delta_n$,
hereafter we only will consider the cases $-\delta_n \le k\le
\delta_1.$ Moreover, the bounds showed in this paper on $a_k(G)$,
for $\delta_n \le k \le \delta_1$, are obtained by supposing that
the graph $G$ contains defensive $k$-alliances.

It was shown in \rm\cite{alliancesOne} that for any graph $G$ of
order  $n$ and minimum degree $\delta_n$,
$$a_{-1}(G)\leq n-\left\lceil\frac{\delta_n}{2}\right\rceil
\quad {\rm and } \quad  a_0(G)\leq
n-\left\lfloor\frac{\delta_n}{2}\right\rfloor .$$
 Here we generalize the previous result
to defensive $k$-alliances  and we obtain lower bounds.

\begin{theorem}\label{cotainfsup}
 For every $k\in \{-\delta_n,\dots,\delta_1\}$,
$$\left\lceil\displaystyle\frac{\delta_n+k+2}{2}\right\rceil \le
a_k(G)\leq
n-\left\lfloor\displaystyle\frac{\delta_n-k}{2}\right\rfloor.$$
 \end{theorem}

\begin{proof}
Let  $X\subseteq V$   be a defensive $k$-alliance in $G$. In this
case, for every $v\in X$ we have
  $$\delta(v) =\delta_X(v)+\delta_{\bar{X}}(v)$$
         $$ \delta(v)  \le \delta_X(v)+ \frac{\delta(v)-k}{2}$$
           $$ \frac{\delta(v)+k}{2} \le \delta_X(v)\le |X|-1$$
         $$ \frac{\delta_n +k+2}{2}\le |X|.$$
Hence, the lower bound follows.

On the other hand, if $X\subseteq V$ is a defensive $k$-alliance in
$G$ for $\delta_n \le k\le \delta_1$, then $a_k(G)\le |X| \le n \le
n-\left\lfloor\frac{\delta_n-k}{2}\right\rfloor.$ Suppose
$-\delta_n\le k \le \delta_n$. Let $S\subseteq V$ be a set of
cardinality $n-\left\lfloor\frac{\delta_n-k}{2}\right\rfloor$. For
every vertex
 $v\in S$ we have $\frac{\delta(v)-k}{2}\ge
\left\lfloor \frac{\delta_n-k}{2} \right\rfloor \ge
\delta_{\bar{S}}(v)$. Hence,
 $S$ is a defensive  $k$-alliance and $a_k(G)\le
|S|=n-\left\lfloor\frac{\delta_n-k}{2}\right\rfloor$.
\end{proof}

We denote by $K_n$ the complete graph of order $n$.

\begin{corollary}\label{coro-kn-ak}
For every $k\in \{1-n, \dots,  n-1\}$, $a_{k}(K_n)=\left\lceil
\displaystyle\frac{n+k+1}{2}\right\rceil$.
\end{corollary}

\begin{theorem}\label{th3}
For every $k,r\in \mathbb{Z}$ such that $ -\delta_n\le k\le
\delta_1$
 and $0 \le r \le \frac{k+\delta_n}{2}$,  $$  a_{_{k-2r}}(G)+r\leq a_{_k}(G).$$
\end{theorem}

\begin{proof}
Let $S\subset V$ be a defensive $k$-alliance of minimum cardinality
in $G$. By Theorem \ref{cotainfsup}, $\frac{\delta_n+k+2}{2}\le
|S|$, then we can take $X\subset S$ such that $|X|=r.$ Hence, for
every $v\in Y=S-X$,
\begin{align*}
\delta_{Y}(v)&=\delta_S(v)-\delta_X(v)\\
                &\ge \delta_{\bar{S}}(v)+k-\delta_X(v) \\
                &= \delta_{\bar{Y}}(v)+k-2\delta_X(v) \\
                &\ge  \delta_{\bar{Y}}(v)+k-2r.
\end{align*}
Therefore, $Y$ is a defensive $(k-2r)$-alliance in $G$ and, as a
consequence, $  a_{k-2r}(G)\leq a_k(G)-r.$
\end{proof}

Notice that, according to the result in Corollary \ref{coro-kn-ak},
the bound for $a_k(G)$ in  Theorem  \ref{th3} is attained for the
complete graph $K_n$ for every $n, k, r$ with its respective
restrictions. From the above theorem we derive some interesting
consequences.

\begin{corollary}\label{coro1} Let $t\in \mathbb{Z}$.
 \begin{itemize}
\item If $\frac{1-\delta_n}{2} \le t \le
\frac{\delta_1-1}{2}$, then  $a_{2t-1}(G)+1\leq a_{2t+1}(G).$
 \item  If $\frac{2-\delta_n}{2} \le t \le
\frac{\delta_1}{2}$, then $ a_{2(t-1)}(G)+1 \leq a_{2t}(G).$
\end{itemize}
\end{corollary}

\begin{corollary}\label{coro2}
For every $k\in \{0,\dots,\delta_n\}$,
\begin{itemize}
\item if $k$ is even, then $a_{-k}(G)+\frac{k}{2}\le a_{0}(G) \leq a_{k}(G)-\frac{k}{2},$
\item if $k$ is odd, then $a_{-k}(G)+\frac{k-1}{2} \le a_{-1}(G)\leq a_{k}(G)-\frac{k+1}{2}.$
\end{itemize}
\end{corollary}

It was shown in \rm\cite{note,alliancesOne} that for any graph $G$
of order  $n$,
\begin{equation}\label{cotasup}
a_{-1}(G)\leq \left\lceil\frac{n}{2}\right\rceil \quad {\rm and
}\quad  a_0(G)\leq \left\lfloor\frac{n}{2}\right\rfloor+1.
\end{equation}

By Corollary \ref{coro2} and (\ref{cotasup}) we obtain the following
result.

\begin{theorem}
For every $k\in \{-\delta_n,\dots ,0\}$, $a_k(G)\le \left\lceil
\displaystyle\frac{n+k+1}{2}\right\rceil.$
 \end{theorem}

Notice that the above bound is attained, for instance, for the
complete graph $G=K_n$.




\section{Algebraic connectivity and defensive $k$-alliance number}
It is well-known that the second smallest Laplacian eigenvalue of a
graph is probably the most important information contained in the
Laplacian spectrum. This eigenvalue, frequently called {\it
algebraic connectivity}, is related to several important graph
invariants and imposes reasonably good bounds on the values of
several parameters of graphs which are very hard to compute.

 The algebraic connectivity of $G$, $\mu$, satisfies the following equality
shown by Fiedler \cite{fiedler} on weighted graphs
\begin{equation}
  \mu=2n \min \left\{ \frac{\sum_{v_i\sim v_j}(w_i-w_j)^2  }
  {\sum_{v_i\in V}\sum_{v_j\in V}(w_i-w_j)^2}: \mbox{\rm $w\neq \alpha{\bf j}$
   for   $\alpha\in \mathbb{R}$ } \right\},     \label{rfiedler}
\end{equation}
where $V=\{v_1, v_2, ..., v_n\}$, ${\bf j}=(1,1,...,1)$ and $w\in
\mathbb{R}^n$.

The following theorem  shows the relationship between the algebraic
connectivity of a graph and its  defensive $k$-alliance number.

\begin{theorem}\label{ak}
For any connected graph $G$  and for every $k\in
\{-\delta_n,\dots,\delta_1\}$,  $$ a_k(G)\ge
\left\lceil\frac{n(\mu+k+1)}{n+\mu}\right\rceil.$$
\end{theorem}

\begin{proof}
If $S$ denotes a defensive $k$-alliance in $G$, then
\begin{equation}\label{cotavecinos}
\delta_{\bar{S}}(v)+k \le |S|-1, \quad \forall v\in S.
\end{equation}
By (\ref{rfiedler}), taking $w\in \mathbb{R}^n$ defined as
$$
w_i= \left\lbrace \begin{array}{ll} 1  & {\rm if }\quad  v_i\in S;
                            \\ 0 &  {\rm  otherwise,} \end{array}
                                                \right  .$$
                                                we have
\begin{equation}\label{fiedlerAlliance}
\mu\le \frac{n\displaystyle\sum_{v\in S} \delta_{\bar{S}}(v)}
{|S|(n-|S|)}.
\end{equation}
Thus,  (\ref{cotavecinos}) and (\ref{fiedlerAlliance}) lead to
\begin{equation}\label{final}
\mu\le \frac{n (|S|-k-1)}{n-|S|}.
\end{equation}
Therefore, solving (\ref{final}) for $|S|$, and considering that it
is an integer, we obtain the bound on $a_k(G)$.\end{proof}

The above bound is sharp as we can see in the following example. As
the algebraic connectivity of the complete  graph $G=K_n$ is
$\mu=n$, the above theorem gives the exact value of
$a_{k}(K_n)=\left\lceil\frac{n+k+1}{2}\right\rceil$.

\begin{theorem}\label{cotafloja}
For any connected graph $G$  and for every $k\in
\{-\delta_n,\dots,\delta_1\}$,
$$a_k(G)\ge \left\lceil\frac{n(\mu-\left\lfloor\frac{\delta_1-k}{2}\right\rfloor)}{\mu}\right\rceil.$$
\end{theorem}

\begin{proof}
If $S$ denotes a defensive $k$-alliance in $G$, then
$\delta_1\ge \delta(v)\ge 2\delta_{\bar{S}}(v)+k, \quad \forall v\in
S.$
Therefore,
\begin{equation} \label{strongGrado}
 \left\lfloor\frac{\delta_1-k}{2}\right\rfloor \ge \delta_{\bar{S}}(v),
\quad \forall v\in S.
\end{equation} Hence, by (\ref{fiedlerAlliance}) and (\ref{strongGrado}) the result follows.\end{proof}

The bound is attained for  every $k$ in  the case of the complete
graph $G=K_n$. 

The reader is referred to \cite{albesiga2} for more details on the
spectral study of offensive alliances and dual alliances.



\section{Defensive $k$-alliance number and line graph}

Hereafter, we denote by ${\cal L}(G)=(V_l,E_l)$ the line graph of a
simple graph $G$. The degree of the vertex $e=\{u,v\}\in V_l$ is
$\delta(e)=\delta(u)+\delta(v)-2$. If the degree sequence of
 $G$ is $\delta_{1}\geq
\delta_{2}\geq \cdots \geq\delta_{n}$, then the maximum degree of
${\cal L}(G)$, denoted by  $\Delta_l$, is bounded by
\begin{equation}
\Delta_l\leq \delta_{1}+ \delta_{2}-2
\end{equation}
and the minimum degree of  ${\cal L}(G)$, denoted by $\delta_l$, is
bounded by
\begin{equation}\label{gradomin}
\delta_l\ge  \delta_{n}+ \delta_{n-1}-2.
\end{equation}

In this section we obtain some results on ${a}_k({\cal L}(G))$ in
terms of the degree sequence of $G$.

\begin{theorem}\label{thcota}
For any simple graph $G$ of maximum degree $\delta_1$, and for every
$k\in \{2(1-\delta_1),\dots,0\}$,
$${a}_k({\cal L}(G))\le
\delta_1+\left\lceil\frac{k}{2}\right\rceil .$$
\end{theorem}

\begin{proof}
Suppose $k$ is even. Let $v\in V(G)$ be a vertex of maximum degree
in $G$ and let  $S_v=\{e\in E: v\in e\}$. Let $Y_k\subset S_v$ such
that $|Y_k|=-\frac{k}{2}$ and let $X_k=S_v-Y_k$. Thus, $\langle
S_v\rangle \cong K_{\delta_1}$ and, as a consequence,
$$\delta_{X_k}(e)=\delta_1-1+\frac{k}{2}\ge
\delta_2-1+\frac{k}{2}\ge \delta_{\bar{X_k}}(e)+k, \quad \forall
e\in X_k.$$ Hence, $X_k\subset V_l$ is a defensive $k$-alliance in
${\cal L}(G)$. So, for $k$  even we have $a_k({\cal L}(G))\le
\delta_1+\left\lceil\frac{k}{2}\right\rceil$. Moreover, if $k$ is
odd, then $a_k({\cal L}(G))\le a_{k+1}({\cal L}(G))\le
\delta_1+\left\lceil\frac{k+1}{2}\right\rceil
=\delta_1+\left\lceil\frac{k}{2}\right\rceil$.
\end{proof}

One advantage of applying the above bound is that it requires only
little information about the graph $G$; just the maximum degree. The
above bound is tight, as we will see below, for any $\delta$-regular
graph, and $k\in \{2(1-\delta),\dots, 0\}$, ${a}_k({\cal L}(G)) =
\delta+\left\lceil\frac{k}{2}\right\rceil .$ Even so, we can improve
this bound for the case of nonregular graphs. The drawback of such
improvement is that we need to know more about $G$.

\begin{theorem}\label{thcota}
Let $G$ be a simple graph, whose degree sequence is $\delta_{1}\geq
\delta_{2}\geq \cdots \geq\delta_{n}$. Let $v\in V$ such that
$\delta(v)=\delta_1$, let $\delta_v=\displaystyle\max_{u: u\sim v}
\{ \delta(u) \}$ and let
$\delta_*=\displaystyle\min_{v:\delta(v)=\delta_1}\{\delta_v\}$. For
every $k\in \{2-\delta_*-\delta_1,\dots,\delta_1-\delta_*\}$,
$${a}_k({\cal L}(G))\le \left\lceil\frac{
\delta_1+\delta_*+k}{2}\right\rceil .$$ Moreover, for every $k\in
\{2-\delta_1-\delta_2,\dots,\delta_1+\delta_2-2\}$,
$$\left\lceil\frac{\delta_{n}+\delta_{n-1}+k}{2}\right\rceil \le
{a}_k({\cal L}(G)).$$
\end{theorem}

\begin{proof}
 Let $v\in V$ a vertex of maximum degree  $\delta(v)=\delta_1$ such that $v$ is adjacent to a
  vertex of degree $\delta_*$. Let
 $S_v=\{e\in E: v\in e\}$. Suppose $\delta_1+\delta_*+k$ is even.
Therefore, taking $S\subset S_v$ such that $|S|=
\frac{\delta_1+\delta_*+k}{2} $, we obtain $\langle S\rangle \cong
K_{ \frac{\delta_1+\delta_*+k}{2} }$. Thus, $\forall e\in S$,
$$\delta_{S}(e)-k= \frac{\delta_1+\delta_*+k}{2} -1-k\ge
\delta_1+\delta_*-1- \frac{\delta_1+\delta_*+k}{2}.$$ On the other
hand, as $\delta(e)=\delta_S(e)+\delta_{\bar{S}}(e)$, we have
$$ \delta_1 + \delta _* - 2\ge
\frac{\delta_1+\delta_*+k}{2}-1+\delta_{\bar{S}}(e)\Leftrightarrow
\delta_1 + \delta _* - 1-\frac{\delta_1+\delta_*+k}{2}\ge
\delta_{\bar{S}}(e).
$$
So, $\delta_{S}(e)-k\ge \delta_{\bar{S}}(e)$ and $S$ is a defensive
$k$-alliance in ${\cal L}(G)$ and, as a consequence, ${a}_{k}({\cal
L}(G))\le \left\lceil\frac{ \delta_1+\delta_*+k}{2}\right\rceil$. If
$\delta_1+\delta_*+k$ is odd, then $\delta_1-\delta_*>k$. Therefore,
${a}_{k}({\cal L}(G))\le {a}_{k+1}({\cal L}(G)) \le
\left\lceil\frac{ \delta_1+\delta_*+k+1}{2}\right\rceil
=\left\lceil\frac{ \delta_1+\delta_*+k}{2}\right\rceil .$

The lower bound follows from Theorem \ref{cotainfsup} and
(\ref{gradomin}).
\end{proof}

\begin{corollary}\label{bonito}
For any $\delta$-regular graph, $\delta>0$, and for every $k\in
\{2(1-\delta),\dots,0\}$,
$$a_k({\cal L}(G))=\delta+\left\lceil\frac{k}{2}\right\rceil .$$
\end{corollary}

We recall that a graph $G=(V, E)$ is a
($\delta_1,\delta_2$-semiregular bipartite graph if the set $V$ can
be partitioned into two disjoint subsets $V_1, V_2$ such that if
$u\sim v$ then $u\in V_1$ and $v\in V_2$ and also,
$\delta(v)=\delta_1$ for every $v\in V_1$ and $\delta(v)=\delta_2$
for every $v\in V_2$.

\begin{corollary}
For any ($\delta_1,\delta_2$)-semiregular bipartite graph $G$,
$\delta_1>\delta_2$, and for every $k\in
\{2-\delta_1-\delta_2,\dots,\delta_1-\delta_2\}$,
 $$a_k({\cal L}(G))=\left\lceil
\frac{\delta_1+\delta_2+k}{2} \right\rceil.$$
\end{corollary}

We should point out that from the results obtained  in the other
sections of this article on $a_k(G)$, we can derive some new results
on $a_k({\cal L}(G))$. The reader is referred to \cite{albesiga8}
for more details on $a_{-1}({\cal L}(G))$ and $a_0({\cal L}(G))$.


\end{document}